\documentclass{amsart}
\usepackage{amsmath}
\usepackage{eurosym}
\usepackage{amsfonts}
\usepackage{amssymb}
\usepackage{amsmath}
\usepackage{amsfonts}
\usepackage{graphicx}
\usepackage{subfig}
\usepackage{eurosym}
\usepackage{amssymb}
\usepackage{amsmath}
\usepackage{amsfonts}
\usepackage{graphicx}
\usepackage{subfig}
\usepackage{float}

\newtheorem{theorem}{Theorem}
\theoremstyle{plain}

\newtheorem{definition}{Definition}

\numberwithin{equation}{section}

\begin{document}
\title[]{The physical approach on the surfaces of rotation in $ E_{2}^{4} $}
\author{Fatma ALMAZ}
\address{faculty of arts and sciences, department of mathematics, batman university, batman /türk\.{ı}ye}
\email{fatma.almaz@batman.edu.tr}
\thanks{This paper is in final form and no version of it will be submitted
for publication elsewhere.}

\begin{abstract}
In this paper, some physical expressions as the specific energy and the specific angular momentum on these surfaces of rotation are investigated with the help of Clairaut's theorem using conditions being geodesic in which the curves can be chosen to be time-like curves, which allows us to constitute the specific energy and specific angular momentum

\end{abstract}

\maketitle

\section{Introduction}

Physical features as energy and momentum that they include the mass as well proportioned factor will instead by changed by the specific features supplied by dividing out the mass. Therefore, since the kinetic energy is $E=m\frac{V^{2}}{2}$, because of feature its motion in space, which the motion is very important in terms of its specific energy and angular momentum in \cite{15,16}. If a force is accountable for this acceleration, that is to say the normal force is perpendicular to the velocity of the particle. Therefore, the specific energy and the speed $V$ must be constant along a geodesic. Because the existence of this constant is a result of the one parameter rotational group of symmetries of the surface, as a constant of the movement introduces a new thing since the surface is invariant under any one parameter group of symmetries, \cite{11}. In \cite{1}, the brief description of rotational surfaces is defined in Galilean 4-space by the authors. In \cite{2}, time-like geodesics are expressed using Clairaut's theorem on the hyperbolic and elliptic rotational surfaces in $E_{2}^{4}$ by the authors. In \cite{3}, the magnetic rotated surfaces are defined in null cone $Q^{2}\subset E_{1}^{3}$ by the authors. In \cite{4}, the 
conditions of being geodesic are expressed on the tube surface using Clairaut's theorem, the specific energy and the angular momentum are defined by the authors. In \cite{5}, different types of rotational surfaces is defined using killing vector field in semi-Euclidean 4-space by the authors. In \cite{8}, A new type of surfaces in Euclidean and Minkowski 4-space is constructed by performing two simultaneous rotations on a planar curve by the authors. Also,
classification theorems of flat double rotational surfaces are proved by the
authors. In \cite{9}, the authors discuss some issues of displaying 2D
surfaces in 4-space, including the behaviour of surface normals under
projection.

\section{Preliminaries}

Let $E_{2}^{4}$ denote the $4-$dimensional pseudo-Euclidean space with
signature $(2,4)$, that is, the real vector space $%
\mathbb{R}
^{4}$ endowed with the metric $\left\langle ,\right\rangle _{E_{2}^{4}}$
which is defined by 
\begin{equation}
\left\langle ,\right\rangle
_{E_{2}^{4}}=-dx_{1}^{2}-dx_{2}^{2}+dx_{3}^{2}+dx_{4}^{2},  \tag{2.1}
\end{equation}%
or 
\begin{equation}
g=%
\begin{bmatrix}
-1 & 0 & 0 & 0 \\ 
0 & -1 & 0 & 0 \\ 
0 & 0 & 1 & 0 \\ 
0 & 0 & 0 & 1%
\end{bmatrix}
\tag{2.2}
\end{equation}%
where $(x_{1},x_{2},x_{3},x_{4})$ is a standard rectangular coordinate
system in $E_{2}^{4}$.

For an arbitrary vector $v\in E_{2}^{4}\backslash \{0\}$ there are one of three characters: it can be space-like if $g(v,v)>0$ or $v=0,$
time-like if $g(v,v)<0$ and null if $g(v,v)=0$ and $v\neq 0.$ Hence, an
arbitrary curve $x(s)$ in $E_{2}^{4}$ can locally be space-like, time-like
or null. Also, the norm of a vector $v$ is given by $\parallel v\parallel =\sqrt{g(v,v)}$ and a space-like or time-like curve $x(s)$ has unit speed, if $g(x^{\prime
},x^{\prime })=\pm 1.$

Let $%
(x_{1},x_{2},x_{3},x_{4}),(y_{1},y_{2},y_{3},y_{4}),(z_{1},z_{2},z_{3},z_{4}) 
$ be any three vectors in $E_{2}^{4}$. The pseudo-Euclidean cross product is
given as 
\begin{equation*}
x\wedge y\wedge z=%
\begin{pmatrix}
-i_{1} & -i_{2} & i_{3} & i_{4} \\ 
x_{1} & x_{2} & x_{3} & x_{4} \\ 
y_{1} & y_{2} & y_{3} & y_{4} \\ 
z_{1} & z_{2} & z_{3} & z_{4}%
\end{pmatrix}%
, \tag{2.3}
\end{equation*}%
where $i_{1}=\left( 1,0,0,0\right) ,i_{2}=\left( 0,1,0,0\right)
,i_{3}=\left( 0,0,1,0\right) ,i_{4}=\left( 0,0,0,1\right) $, \cite{7,9,12,14}.

\begin{definition}
Let $W$ be a vector field on a smooth manifold $M$ and $\psi_{t}$ be the
local flow generated by $W$. For each $t\in%
\mathbb{R}
,$ the map $\psi_{t}$ is diffeomorphism of $M$ and given a function $f$ on $%
M $, one considers the Pull-back $\psi_{t}f$, the Lie derivative of the
function $f$ as defined as to $W$ by%
\begin{equation}
L_{_{W}}f=\underset{t\longrightarrow0}{\lim}\underset{}{\left( \frac{\psi
_{t}f-f}{t}\right) =\frac{d\psi_{t}f}{dt}_{t=0}}.  \tag{2.4}
\end{equation}

Let $g_{\xi\varrho}$ be any pseudo-Riemann metric, then the derivative is
given as 
\begin{equation*}
L_{_{W}}g_{\xi\varrho}=g_{\xi\varrho,z}W^{z}+g_{\xi
z}W_{,\varrho}^{z}+g_{z\varrho}W_{,\xi}^{z}.
\end{equation*}

In Cartesian coordinates in Euclidean spaces where $g_{\xi\varrho,z}=0,$ and
the Lie derivative is given by%
\begin{equation*}
L_{_{W}}g_{\xi\varrho}=g_{\xi z}W_{,\varrho}^{z}+g_{z\varrho}W_{,\xi}^{z},
\end{equation*}
the vector $W$ generates a Killing field if and if only 
\begin{equation*}
L_{_{W}}g=0.
\end{equation*}
\end{definition}

in \cite{6,10,11,17}.


\begin{theorem}
Let the pseudo-Euclidean group be a subgroup of the diffeomorphisms group in 
$E_{2}^{4}$ and let $W$ be vector field which generate the isometries. Then,
the killing vector field associated with the metric $g$ is given as%
\begin{equation*}
W(\xi ,\varrho ,\vartheta ,\eta )=a\left( \eta \partial \xi +\xi \partial
\eta \right) +b\left( \vartheta \partial \varrho +\varrho \partial \vartheta
\right) +c\left( \vartheta \partial \xi +\xi \partial \vartheta \right)
\end{equation*}%
\begin{equation*}
+d(\eta \partial \varrho +\varrho \partial \eta )+e(\vartheta \partial \eta
-\eta \partial \vartheta )+f\left( \xi \partial \varrho -\varrho \partial
\xi \right) ,
\end{equation*}%
where $a,b,c,d,e,f\in 
\mathbb{R}
_{0}^{+},$ \cite{5}.
\end{theorem}

\begin{theorem}
Let $W(\xi ,\varrho ,\vartheta ,\eta )$ be the killing vector field and let $%
\gamma=(f_{1},f_{2},f_{3},f_{4})$ be a curve in $E_{2}^{4}$, then the
surfaces of rotation are given as follows

\begin{enumerate}
\item For the rotations $\Omega _{1}=\vartheta \partial \xi +\xi \partial
\vartheta $ and $\Omega _{4}=\eta \partial \varrho +\varrho \partial \eta ,$
the hyperbolic surface of rotation is given as%
\begin{equation*}
S_{14}(x,\alpha ,s)=\left( 
\begin{array}{c}
f_{1}\cosh x+f_{3}\sinh x,f_{2}\cosh \alpha +f_{4}\sinh \alpha , \\ 
f_{1}\sinh x+f_{3}\cosh x,f_{2}\sinh \alpha +f_{4}\cosh \alpha 
\end{array}%
\right) 
\end{equation*}%
and for the planar curve $\gamma (s)=(f_{1}(s),0,0,f_{4}(s))$ the Gaussian
curvature $K$ and the mean curvature vector $H$ of the rotational surface $%
S_{14}(x(t),\alpha (t),s)=\left( f_{1}\cosh x,f_{4}\sinh \alpha ,f_{1}\sinh
x,f_{4}\cosh \alpha \right) $ are given as%
\begin{equation*}
K=\frac{\left( f_{1}^{\prime }f_{4}^{{}}-f_{1}f_{4}^{\prime }\right)
^{2}\left( \overset{.}{x}\overset{.}{\alpha }\right) ^{2}}{f_{4}^{2}\overset{%
.}{\alpha }^{2}-f_{1}^{2}\overset{.}{x}^{2}}+\frac{\left( f_{1}^{\prime
}f_{4}\overset{.}{\alpha }^{2}-f_{4}^{\prime }f_{1}\overset{.}{x}^{2}\right)
\left( f_{1}^{\prime }f_{4}^{\prime \prime }-f_{1}^{\prime \prime
}f_{4}^{\prime }\right) }{-f_{1}^{\prime 2}+f_{4}^{\prime 2}},
\end{equation*}%
\begin{equation*}
H=\{\frac{f_{1}f_{4}\left( \overset{..}{x}\overset{.}{\alpha }+\overset{.}{x}%
\overset{..}{\alpha }\right) }{2\sqrt{f_{4}^{2}\overset{.}{\alpha }%
^{2}-f_{1}^{2}\overset{.}{x}^{2}}}+\frac{f_{4}^{\prime }f_{1}\overset{.}{x}%
^{2}-f_{1}^{\prime }f_{4}\overset{.}{\alpha }^{2}}{2\sqrt{-f_{1}^{\prime
2}+f_{4}^{\prime 2}}}\}e_{3}+\frac{\left( f_{1}^{\prime }f_{4}^{\prime
\prime }-f_{1}^{\prime \prime }f_{4}^{\prime }\right) }{2\sqrt{%
-f_{1}^{\prime 2}+f_{4}^{\prime 2}}}e_{4}
\end{equation*}%
where\\
$e_{3}=\frac{\left( f_{4}\overset{.}{\alpha }\sinh x,f_{1}\overset{.}{x%
}\cosh \alpha ,f_{4}\overset{.}{\alpha }\cosh x,f_{1}\overset{.}{x}\sinh
\alpha \right) }{\sqrt{f_{4}^{2}\overset{.}{\alpha }^{2}-f_{1}^{2}\overset{.}%
{x}^{2}}}$,
$e_{4}=\frac{\left( f_{4}^{\prime }\cosh x,f_{1}^{\prime }\sinh
\alpha ,f_{4}^{\prime }\sinh x,f_{1}^{\prime }\cosh \alpha \right) }{\sqrt{%
-f_{1}^{\prime 2}+f_{4}^{\prime 2}}}.$

\item For the rotations $\Omega _{2}=\eta \partial \xi +\xi \partial \eta $
and $\Omega _{3}=\vartheta \partial \varrho +\varrho \partial \vartheta ,$
the hyperbolic surface of rotation is given as 
\begin{equation*}
S_{23}(y,z,s)=\left( 
\begin{array}{c}
f_{1}\cosh y+f_{4}\sinh y,f_{2}\cosh z+f_{3}\sinh z, \\ 
f_{2}\sinh z+f_{3}\cosh z,f_{1}\sinh y+f_{4}\cosh y%
\end{array}%
\right) .
\end{equation*}%
and for the planar curve $\gamma (s)=(f_{1}(s),f_{2}(s),0,0)$ the Gaussian
curvature $K$ and the mean curvature vector $H$ of the rotational surface $%
S_{23}(y(t),z(t),s)=\left( f_{1}\cosh y,f_{2}\cosh z,f_{2}\sinh z,f_{1}\sinh
y\right) $ are given as%
\begin{equation*}
K=-\left( 
\begin{array}{c}
\frac{\left( f_{1}f_{2}^{\prime }+f_{1}^{\prime }f_{2}^{{}}\right)
^{2}\left( \overset{.}{y}\overset{.}{z}\right) ^{2}}{f_{2}^{2}\overset{.}{z}^{2}%
+f_{1}^{2}\overset{.}{y}^{2}}+ \\ 
\frac{\left( f_{1}^{{}}f_{2}^{\prime }\overset{.}{y}^{2}+f_{1}^{\prime }f_{2}%
\overset{.}{z}^{2}\right) \left( f_{1}^{\prime \prime }f_{2}^{\prime
}+f_{1}^{\prime }f_{2}^{\prime \prime }\right) }{f_{1}^{\prime
2}+f_{2}^{\prime 2}}%
\end{array}%
\right) ;H=\left( 
\begin{array}{c}
\frac{f_{1}^{{}}f_{2}(\overset{.}{y}\overset{..}{z}+\overset{..}{y}\overset{.%
}{z})}{2\sqrt{f_{2}^{2}\overset{.}{z}^{2}+f_{1}^{2}\overset{.}{y}^{2}}}e_{3} \\ 
+\frac{f_{1}^{{}}f_{2}^{\prime }\overset{.}{y}^{2}+f_{1}^{\prime }f_{2}%
\overset{.}{z}^{2}-f_{1}^{\prime \prime }f_{2}^{\prime }-f_{1}^{\prime
}f_{2}^{\prime \prime }}{2\sqrt{f_{1}^{\prime 2}+f_{2}^{\prime 2}}}e_{4}%
\end{array}%
\right) ,
\end{equation*}%
where\\
$e_{3}=\frac{\left( f_{2}\overset{.}{z}\sinh y,f_{1}\overset{.}{y}%
\sinh z,f_{1}\overset{.}{y}\cosh z,f_{2}\overset{.}{z}\cosh y\right) }{\sqrt{%
f_{2}^{2}\overset{.}{z}^{2}+f_{1}^{2}\overset{.}{y}^{2}}}$,
$e_{4}=\frac{\left(
f_{2}^{\prime }\cosh y,f_{1}^{\prime }\cosh z,f_{1}^{\prime }\sinh
z,f_{2}^{\prime }\sinh y\right) }{\sqrt{f_{1}^{\prime 2}+f_{2}^{\prime 2}}}$

\item For the rotations $\Omega _{5}=\xi \partial \varrho -\varrho \partial
\xi $ and $\Omega _{6}=\vartheta \partial \eta -\eta \partial \vartheta ,$
the elliptic surface of rotation is given as 
\begin{equation*}
S_{56}(\beta ,\theta ,s)=\left( 
\begin{array}{c}
f_{1}\cos \beta +f_{2}\sin \beta ,-f_{1}\sin \beta +f_{2}\cos \beta , \\ 
f_{3}\cos \theta +f_{4}\sin \theta ,-f_{3}\sin \theta +f_{4}\cos \theta 
\end{array}%
\right) ,
\end{equation*}%
and for the planar curve $\gamma (s)=(0,f_{2}(s),0,f_{4}(s))$ the Gaussian
curvature $K$ and the mean curvature vector $H$ of the rotational surface $%
S_{56}(\beta \left( t\right) ,\theta \left( t\right) ,s)=\left( f_{2}\sin
\beta ,f_{2}\cos \beta ,f_{4}\sin \theta ,f_{4}\cos \theta \right) $ are
given as 
\begin{equation*}
K=-\left( 
\begin{array}{c}
\frac{\left( f_{2}^{\prime }f_{4}-f_{2}f_{4}^{\prime }\right) ^{2}\left( 
\overset{.}{\beta }\overset{.}{\theta }\right) ^{2}}{-f_{2}^{2}\overset{.}{%
\beta }^{2}+f_{4}^{2}\overset{.}{\theta }^{2}}+
\frac{\left( -f_{2}^{\prime \prime }f_{4}^{\prime }+f_{2}^{\prime
}f_{4}^{\prime \prime }\right) (f_{4}^{\prime }f_{2}\overset{.}{\beta }%
^{2}-f_{2}^{\prime }f_{4}\overset{.}{\theta }^{2})^{2}}{-f_{2}^{\prime
2}+f_{4}^{\prime 2}}%
\end{array}%
\right) ;
\end{equation*}
\begin{equation*}
H=\frac{f_{4}f_{2}\left( \overset{.}{\beta }\overset{..}{\theta }-\overset{.}%
{\theta }\overset{..}{\beta }\right) }{2\sqrt{f_{4}^{2}\overset{.}{\theta }^{2}%
-f_{2}^{2}\overset{.}{\beta }^{2}}}e_{3}+\frac{\left( f_{4}^{\prime }f_{2}%
\overset{.}{\beta }^{2}-f_{2}^{\prime }f_{4}\overset{.}{\theta }%
^{2}+f_{2}^{\prime \prime }f_{4}^{\prime }-f_{2}^{\prime }f_{4}^{\prime
\prime }\right) }{2\sqrt{f_{4}^{\prime 2}-f_{2}^{\prime 2}}}e_{4}
\end{equation*}
where\\
$e_{3}=\frac{\left( -f_{4}\overset{.}{\theta }\cos \beta ,f_{4}\overset%
{.}{\theta }\sin \beta ,-f_{2}\overset{.}{\beta }\cos \theta ,f_{2}\overset{.%
}{\beta }\sin \theta \right) }{\sqrt{-f_{2}^{2}\overset{.}{\beta }%
^{2}+f_{4}^{2}\overset{.}{\theta }^{2}}}$,
$e_{4}=\frac{\left( f_{4}^{\prime }\sin\beta ,f_{4}^{\prime }\cos \beta ,f_{2}^{\prime }\sin \theta ,f_{2}^{\prime
}\cos \theta \right) }{\sqrt{-f_{2}^{\prime 2}+f_{4}^{\prime 2}}}$; $ -\infty
<x,y,z,\alpha ,\beta ,\theta <\infty ,s\in I$ and $f_{i}\in C^{\infty }$, \cite{5}.
\end{enumerate}
\end{theorem}

\begin{theorem}
Let $\gamma (t)=(f_{1}(t),0,0,f_{4}(t))$(or $\gamma
(t)=(0,f_{2}(t),f_{3}(t),0)$)$,f_{i}\in C^{\infty }$ be a time-like geodesic
curve on the hyperbolic surface of rotation $S_{14}$ in the $%
E_{2}^{4} $, let $f_{1}$ and $f_{4}$ be the distance functions from the axis
of rotation to a point on the surface. Therefore, $2f_{1}\cos \varphi _{1}$
and $-2f_{4}\cosh \theta _{1}\sin \varphi _{1}$ are constant along the curve 
$\gamma $ where $\varphi _{1}$ and $\theta _{1}$ are the angles between the
meridians of the surface and the time-like geodesic $\gamma $. Conversely,
if $2f_{1}\cos \varphi _{1}$ and $-2f_{4}\cosh \theta _{1}\sin \varphi _{1}$
are constant along $\gamma $, if no part of some parallels of the surface of
rotation, then $\gamma $ is time-like geodesic \cite{2}.
\end{theorem}


\begin{theorem}
\cite{2}, The general equation of geodesics on the rotational surface $%
S_{14}\subset E_{2}^{4}$, and for the parameters $\overset{.}{x}=\frac{%
1}{f_{1}}\cos \varphi_{1}$ and $\overset{.}{\alpha }=\frac{1}{f_{4}}\cosh
\theta_{1} \sin\varphi_{1}$, are given by 
\begin{equation*}
\frac{dt}{dx}^{{}}=f_{1}\sqrt{1-\cosh ^{2}\theta_{1} \tan ^{2}\varphi_{1}
-L\sec^{2}\varphi_{1} }
\end{equation*}
or 
\begin{equation*}
\frac{dt}{d\alpha }^{{}}=f_{2}\sqrt{\cot ^{2} \varphi_{1} \tan h^{2}
\theta_{1} -L \ sec h^{2} \varphi_{1} \ cosec^{2} \varphi_{1}}.
\end{equation*}
\end{theorem}

\begin{theorem}
Let $\gamma (t)=(f_{1}(t),f_{2}(t),0,0)$(or $\gamma
(t)=(0,0,f_{3}(t),f_{4}(t))$)$,f_{i}\in C^{\infty }$ be a time-like geodesic
curve on the hyperbolic surface of rotation $S_{23}$ in the $%
E_{2}^{4} $, and let $f_{1}$ and $f_{2}$ be the distance functions from the
axis of rotation to a point on the surface. Then, $2f_{1}\cos \theta
_{2}\sinh \varphi _{2}$ and $2f_{2}\sin \theta _{2}\sinh \varphi _{2}$ are
constant along the curve $\gamma $ where $\varphi _{2}$ and $\theta _{2}$
are the angles between the meridians of the surface and the time-like
geodesic curve $\gamma $. Conversely, if $2f_{1}\cos \theta _{2}\sinh
\varphi _{2}$ and $2f_{2}\sin \theta _{2}\sinh \varphi _{2}$ are constant
along the curve $\gamma $, if no part of some parallels of the surface of
rotation, then $\gamma $ is time-like geodesic \cite{2}.
\end{theorem}


\begin{theorem}
\cite{2}, The general equation of geodesics on the rotational surface $%
S_{23}\subset E_{2}^{4} $, and for the parameters $\overset{.}{y}=%
\frac{\cos \theta_{2} \sinh \varphi_{2} }{f_{1}}$ and $\overset{.}{z}=\frac{%
\sinh \varphi_{2} \sin \theta_{2} }{f_{2}}$, are given by 
\begin{equation*}
\frac{dt}{dx}^{{}}=\frac{f_{1}}{\cos \theta_{2} \sinh \varphi_{2} }\sqrt{%
\sinh ^{2}\varphi_{2} -L};\frac{dt}{dz}^{{}}=\frac{f_{2}}{\sinh \varphi_{2}
\sin \theta_{2} }\sqrt{\sinh ^{2}\varphi_{2} -L}.
\end{equation*}
\end{theorem}

\begin{figure}[htp]
\centering
\subfloat[ \label{fig14ax}] {\
\includegraphics[width=7cm,height=7.5cm]{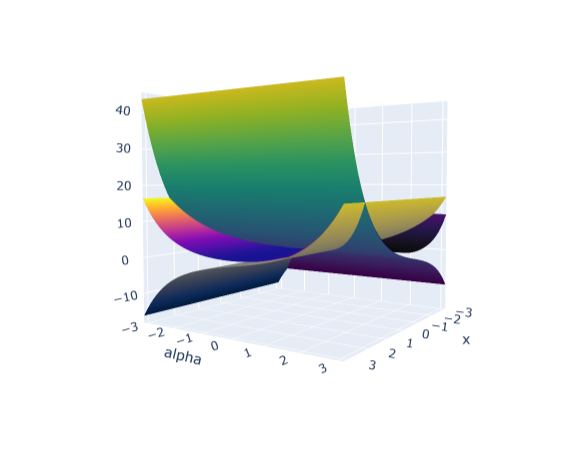} } \hspace*{.1cm}
\subfloat[  \label{fig23a}] {\
\includegraphics[width=6cm,height=7.5cm]{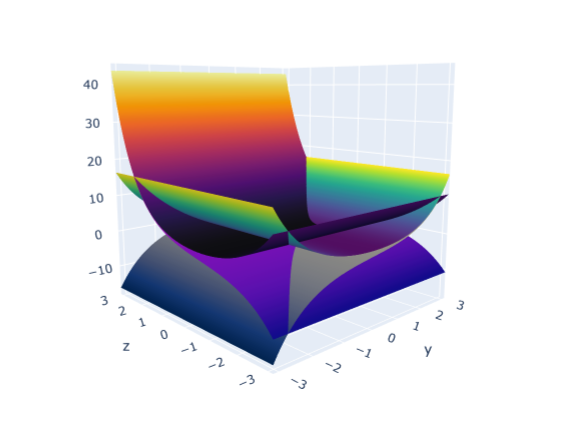}} \newline\caption{{Graphics of
hyperbolic rotational surfaces $ S_{14}(x,\alpha,s) $ and $ S_{23}(y,z,s) $ generated by the curve $\gamma(s)=(2coss,2sins,3s,0)$} }%
\label{Fig1}%
\subfloat[ \label{fig56a}] {\
\includegraphics[width=13cm,height=8.5cm]{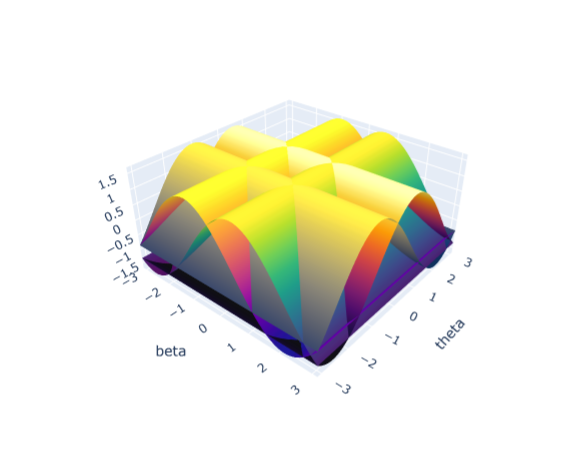} } \hspace*{.1cm}
\newline\caption{{Graphic of
elliptic rotational surface $ S_{56}(\beta,\theta,s) $ generated by the curve $\gamma(s)=(2coshs,2sinhs,2coshs,2sinhs)$} }%
\label{Fig2}%
\end{figure}

\begin{theorem}
Let $\gamma (t)=(0,f_{2}(t),0,f_{4}(t))$(or $\gamma
(t)=(f_{1}(t),0,f_{3}(t),0)$)$,f_{i}\in C^{\infty }$ be a time-like geodesic
curve on the elliptic surface of rotation $S_{56}\subset E_{2}^{4}$,
and let $f_{2}$ and $f_{4}$ be the distance functions from the axis of
rotation to a point on the surface. Then, $2f_{2}\sin \varphi _{3}\cosh
\theta _{3}$ and $2f_{4}\sinh \theta _{3}\sin \varphi _{3}$ are constant
along the curve $\gamma $ where $\varphi _{3}$ and $\theta _{3}$ are the
angles between the meridians of the surface and the time-like geodesic curve 
$\gamma $. Conversely, if $2f_{2}\sin \varphi _{3}\cosh \theta _{3}$ and $%
2f_{4}\sinh \theta _{3}\sin \varphi _{3}$ are constant along the curve $%
\gamma $, if no part of some parallels of the surface of rotation, then $%
\gamma $ is time-like geodesic curve \cite{2}.
\end{theorem}


\begin{theorem}
\cite{2}, The general equation of geodesics on the rotational surface $%
S_{56}\subset E_{2}^{4} $, and for the parameters $\overset{.}{\beta }=%
\frac{\sin \varphi _{3}\cosh \theta _{3}}{f_{2}}$ and $\overset{.}{\upsilon }%
=\frac{\sinh \theta_{3}\sin \varphi _{3}}{f_{4}} $, are given by 
\begin{equation*}
\frac{dt}{d\beta }^{{}}=i\frac{f_{2}\sqrt{L+\sin ^{2}\varphi _{3}}}{\sin
\varphi _{3}\cosh \theta _{3}};\frac{dt}{d\upsilon }^{{}}=i\frac{f_{4}}{%
\sinh \theta _{3}\sin \varphi _{3}}\sqrt{\sin ^{2}\varphi _{3}+L}.
\end{equation*}
\end{theorem}

\section{Physical approach on the surfaces of rotation in $E_{2}^{4}$}

In this section, 
by using the variational approach, which produces the geodesics by
extremizing an action functional on the space of all curves connecting any
two fixed points on the surfaces of rotation ( the hyperbolic surfaces of
rotation $S_{14}=\Upsilon ^{1}(x,\alpha ,t)$, $S_{23}=\Upsilon ^{2}(y,z,t)$ and the elliptic surface of rotation $S_{56}=\Upsilon ^{3}(\beta ,\theta ,t)$). Hence, one can go a step further than all the Riemann geometry discussions about covariant differentiation and parallel transport.

1) For the hyperbolic surface of rotation $\Upsilon^{1}$; one will try to
obtain specific energy equations on this surface. Then, let $%
\Upsilon^{1}(x\left( s\right),\alpha \left( s\right) ,t(s))$ be a
parametrized curve on the surface, which is the integral of the length of the tangent vector in any parametrization of the curve, the speed function is given as 
\begin{equation*}
I^{1}_{1}=\int ds=\int \frac{ds}{d\pi }d\pi =\int \sqrt{\left( \frac{dx}{%
d\pi }\right) ^{2}+\left( \frac{d\alpha }{d\pi }\right) ^{2}+\left( \frac{dt%
}{d\pi }\right) ^{2}}d\pi,
\end{equation*}%
and the integral of half the length squared of the tangent vector%
\begin{equation*}
I^{1}_{2}=\frac{1}{2}\int \left( \frac{ds}{d\pi }\right) ^{2}d\pi =\frac{1}{2%
}\int \left( \left( f_{1}\frac{dx}{d\pi }\right) ^{2}-\left( f_{4}\frac{%
d\alpha }{d\pi }\right) ^{2}-\left( \frac{dt}{d\pi }\right) ^{2}\right) d\pi,
\end{equation*}%
which the speed $\frac{ds}{d\pi }$ is constant and integrate is said to be a Lagrangian function. The second Lagrangian function is 
\begin{equation*}
L^{1}_{2}=\left( x,\alpha ,t,\frac{dx}{d\pi },\frac{d\alpha }{d\pi },\frac{dt%
}{d\pi }\right) =\frac{1}{2}\left( f_{1}\overset{.}{x}\right) ^{2}-\frac{1}{2%
}\left( f_{4}\overset{.}{\alpha }\right) ^{2}-\frac{1}{2}\left( \overset{.}{t%
}\right) ^{2}=E^{1},
\end{equation*}
which is the energy function, while the first Lagrangian $L^{1}_{1}=\frac{ds}{d\pi }$ is speed function
given as 
\begin{equation*}
L^{1}_{1}=\left( x,\alpha ,t,\overset{.}{x},\overset{.}{\alpha },\overset{.}{%
t}\right) =\sqrt{\left( f_{1}\overset{.}{x}\right) ^{2}-\left( f_{4}\overset{%
.}{\alpha }\right) ^{2}-\left( \overset{.}{t}\right) ^{2}}.
\end{equation*}

Both are independent of the azimuthal angle because of the rotational
invariance of the problem. Then, the Lagrange equation of motion of a
particle, analogues to equation defined in terms of the Lagrangian $%
L^{i}_{2} $ with the non-scalar time variable $\pi$ as the parameter, are
given by 
\begin{equation*}
\frac{d}{d\pi }(\frac{\partial L^{i}_{2}}{\partial \left( \frac{\partial
a^{j}}{\partial \pi }\right) })=\frac{\partial L^{i}_{2}}{\partial a^{j}}%
;i,j=1,2,3
\end{equation*}%
\cite{16}, with the angular equation giving the constancy of the angular momentum $l_{i}=\frac{\partial L^{i}_{2}}{\partial 
\overset{.}{a}^{j}}$. Hence, the constancy of the momentum conjugate to $a$
is written as 
\begin{equation*}
p^{a}=\frac{\partial L^{i}_{2}}{\partial \left( \frac{\partial a^{j}}{%
\partial \pi }\right) }
\end{equation*}
\cite{16}, and let us now calculate the total time derivative of the
Lagrangian $L^{i}_{2}$ as follows%
\begin{equation*}
\frac{\partial L^{i}_{2}}{\partial \pi }=\frac{\partial L^{i}_{2}}{\partial
a^{j}}\frac{\partial a^{j}}{\partial \pi }+\frac{\partial L^{i}_{2}}{%
\partial v^{j}}\frac{\partial v^{j}}{\partial \pi };\frac{\partial a^{j}}{%
\partial \pi }=v^{j},
\end{equation*}%
by using the equations of motion and the definition of the three dimensional
velocity can be written as 
\begin{equation*}
\frac{\partial }{\partial \pi }\left( \frac{\partial L^{i}_{2}}{\partial
v^{j}}v^{j}-L^{i}_{2}\right) =0;\frac{\partial L^{i}_{2}}{\partial v^{j}}%
v^{j}-L^{i}_{2}=constant.
\end{equation*}

For the curve $\Upsilon ^{1}(x\left( s\right) ,\alpha \left( s\right) ,t(s))$,  the tangent vector of this curve can be evaluated by using the chain rule and theorem 4, one gets 
\begin{equation}
\frac{d\Upsilon ^{1}(x\left( s\right) ,\alpha \left( s\right) ,t(s))}{ds}=%
\frac{dx\left( s\right) }{ds}\Upsilon _{x}^{1}+\frac{d\alpha \left( s\right) 
}{ds}\Upsilon _{\alpha}^{1}+\frac{dt\left( s\right) }{ds}\Upsilon _{t}^{1}; 
\tag{3.1}
\end{equation}%
\begin{equation}
\overset{.}{\gamma }=N_{x}\cos \varphi _{1}+N_{x}^{\bot }\sin \varphi _{1}=%
\overset{.}{x}\Upsilon _{x}^{1}+\overset{.}{\alpha }\Upsilon _{\alpha }^{1}+%
\overset{.}{t}\Upsilon _{t}^{1}  \tag{3.2}
\end{equation}%
\begin{equation*}
\overset{.}{\gamma }=f_{1}N_{x}\overset{.}{x}+\left( f_{2}N_{\alpha }\overset%
{.}{\alpha }+\overset{.}{t}N_{t}\right) =N_{x}\cos \varphi _{1}+N_{x}^{\bot
}\sin \varphi _{1};
\end{equation*}%
\begin{equation}
=\cos \varphi _{1}N_{x}+\cosh \theta _{1}\sin \varphi _{1}N_{\alpha }+\sinh
\theta _{1}\sin \varphi _{1}N_{t}.  \tag{3.3}
\end{equation}

The tangent vector of the geodesic is given as 
\begin{equation*}
\overset{\rightarrow }{V_{1}}=\frac{d\Upsilon ^{1}}{ds}=V_{1}^{x}\Upsilon
_{x}^{1}+V_{1}^{\alpha }\Upsilon _{\alpha }^{1}+V_{1}^{t}\Upsilon _{t}^{1}
\end{equation*}%
and one can write component vectors notation for components
with respect to the basis vectors $\Upsilon _{x}^{1},\Upsilon _{\alpha
}^{1},\Upsilon _{t}^{1}$ as 
\begin{equation*}
V_{j}^{i}=\frac{dz^{j}}{ds};\left\langle V_{j}^{x},V_{j}^{\alpha
}\right\rangle =\left\langle \frac{dx}{ds},\frac{d\alpha }{ds}\right\rangle
\end{equation*}%
and $V_{1}=\left\langle \overrightarrow{V_{1}},\overrightarrow{%
V_{1}}\right\rangle $ $^{1/2}=\sqrt{g_{ij}\frac{dz^{i}}{ds}\frac{dz^{j}}{ds}}
$ is the speed, which is just the time rate of change of the arc length
along the curve $\gamma $.

Think that $V_{1}^{x^{\ast }}=f_{1}V_{1}^{x}$ $=V_{1}\cos \varphi_{1} $ and $%
V_{1}^{\alpha ^{\ast }}=f_{4}V_{1}^{\alpha }$ $=V_{1}\cosh \theta_{1} \sin
\varphi_{1} $ are just the radial vertical velocity while $V_{1}^{t}$ is the
horizontal angular velocity and $V_{1}^{t^{\ast }}=V_{1}^{t}$ $=V_{1}\sinh
\theta_{1} \sin \varphi_{1} $ is horizontal component of the velocity
vector. The velocity can be represented according to polar coordinates in the
tangent plane to make explicit its magnitude and slope angle with respect to
the radial direction on the surface.

One represents the orthonormal components in terms of the usual polar
coordinate variables in this velocity plane in which $V_{1}^{x^{\ast }}$ is
along the first axis, $V_{1}^{\alpha ^{\ast }}$ is along the second axis and 
$V_{1}^{t^{\ast }}$ is along the third axis.

The speed plays the role of the radial variable in this velocity plane,
while the angles $\theta_{1} $ and $\varphi_{1} $ give the direction of the
velocity according to the direction $\Upsilon _{x^{\ast }}^{1}$ in the
counter clockwise sense in this plane. Also, one can say that the speed is
constant along the geodesic.

It is to understand the system of two second order geodesic equations that
one can use a standard physics technique of partially integrating them and
so lessen them to two first order equations by using two constants of the
movement. From the mass $m$ of the point particle is insufficient in this study. Thus, the specific
kinetic energy can be written given as follows 
\begin{equation*}
E^{1}_{_{\substack{ specific  \\ energy}}}=\frac{1}{2}V_{1}^{2}=\frac{1}{2}%
\left( V_{1}^{2}\cos ^{2}\varphi_{1} +V_{1}^{2}\cosh ^{2}\theta_{1} \sin
^{2}\varphi_{1} -V_{1}^{2}\sinh ^{2}\theta_{1} \sin ^{2}\varphi_{1} \right)
\end{equation*}%
\begin{equation}
=\frac{1}{2}\left( f_{1}\frac{dx}{ds}\right) ^{2}-\frac{1}{2}\left( f_{4}%
\frac{d\alpha }{ds}\right) ^{2}-\frac{1}{2}\left( \frac{dt}{ds}\right) ^{2},
\tag{3.4}
\end{equation}%
then in the physics approach the specific energy and speed are constant along a geodesic. Therefore, specific kinetic energy $E^{1}$ and $V_{1}=\sqrt{2E^{1}}$ 
must be constant along a geodesic.

From Theorem 4 and Theorem 5, for $x=\int \frac{1}{f_{1}}\cos \varphi_{1} ds$
and $\alpha =\int \frac{1}{f_{4}}\cosh \theta_{1} \sin \varphi_{1} ds$ one
can write exactly as in the case of circular motion around an axis with
radius 
\begin{equation*}
\left\Vert \overset{\rightarrow }{R_{1}}\right\Vert =f_{1}\text{ and }%
\left\Vert \overset{\rightarrow }{R_{2}}\right\Vert =f_{4}\text{ or }\overset%
{\rightarrow }{R_{1}}=f_{1}\overset{\rightarrow }{e_{1}}\text{ and }\overset{%
\rightarrow }{R_{2}}=f_{4}\overset{\rightarrow }{e_{2}}.
\end{equation*}

That is, to know the velocity $V_{1}^{x^{\ast }}$ $=V_{1}\cos \varphi_{1}
=f_{1}\frac{dx}{ds}$ and the velocity $V_{1}^{\alpha ^{\ast }}$ $%
=-V_{1}\cosh \theta_{1} \sin \varphi_{1} =f_{4}\frac{d\alpha }{ds}$, the
velocity $V_{1}^{t^{\ast }}=V_{1}\sinh \theta_{1} \sin \varphi_{1}$ $=\frac{%
dt}{ds}$ in the angular direction multiplied by the radius $f_{2}$ and $%
f_{4}.$ Physically, since the second geodesic equation, one writes the
following equations 
\begin{equation*}
l_{_{\substack{ specific\text{ }angular  \\ momentum}}}=\frac{\partial
L^{1}_{2}}{\partial \overset{.}{t}}=-2\overset{.}{t}=-2\sinh \theta_{1} \sin
\varphi_{1} V_{1}\Rightarrow =\frac{-l_{1}}{2}=\overset{.}{t}.
\end{equation*}

The specific angular momentum about the axis of symmetry is constant along a
geodesic. This expression can be used to rewrite the variable angular
velocity $dx/ds$ and $d\alpha /ds$ in the specific energy formula, to obtain the constant specific kinetic energy that is given according to the radial motion and another constant of the motion is given as 
\begin{equation}
E_{_{\substack{ specific  \\ energy}}}=\frac{V_{1}^{2}}{2}\left( \cos
^{2}\varphi_{1} -\cosh ^{2}\theta_{1} \sin ^{2}\varphi_{1} \right) -\frac{%
l_{1}^{2}}{8}.  \tag{3.5}
\end{equation}

2) For the hyperbolic surface of rotation $\Upsilon ^{2}(y\left( s\right)
,z\left( s\right) ,t(s));$ similarly, one can write the speed function 
\begin{equation*}
I_{1}^{2}=\int ds=\int \frac{ds}{d\pi }d\pi =\int \sqrt{\left( \frac{dy}{%
d\pi }\right) ^{2}+\left( \frac{dz}{d\pi }\right) ^{2}+\left( \frac{dt}{d\pi 
}\right) ^{2}}d\pi ,
\end{equation*}%
which is clearly independent of a change of parametrization or the integral
of half the length squared of the tangent vector%
\begin{equation*}
I_{2}^{2}=\frac{1}{2}\int \left( \frac{ds}{d\pi }\right) ^{2}d\pi =\frac{1}{2%
}\int \left( \left( f_{1}\frac{dy}{d\pi }\right) ^{2}+\left( f_{2}\frac{dz}{%
d\pi }\right) ^{2}-\left( \frac{dt}{d\pi }\right) ^{2}\right) d\pi ,
\end{equation*}%
which is equivalent to the previous case only for affinely parametrized
curves for the speed $\frac{ds}{d\pi }$ being constant and is given as 
\begin{equation*}
L_{1}^{2}=\left( y,z,t,\overset{.}{y},\overset{.}{z},\overset{.}{t}\right) =%
\sqrt{\left( f_{1}\overset{.}{y}\right) ^{2}+\left( f_{2}\overset{.}{z}%
\right) ^{2}-\left( \overset{.}{t}\right) ^{2}}
\end{equation*}%
and the integrate is a Lagrangian function that is a function of the curve
and its tangent vector. The second Lagrangian function is the energy
function given as 
\begin{equation*}
L_{2}^{2}=\left( y,z,t,\frac{dy}{d\pi },\frac{dz}{d\pi },\frac{dt}{d\pi }%
\right) =\frac{1}{2}\left( f_{1}\overset{.}{y}\right) ^{2}+\frac{1}{2}\left(
f_{2}\overset{.}{z}\right) ^{2}-\frac{1}{2}\left( \overset{.}{t}\right)
^{2}=E^{2}.
\end{equation*}

Also, in order to calculate the derivative of this tangent vector along the
curve $\Upsilon ^{2}(y\left( s\right) ,z\left( s\right) ,t(s))$. Thus, the
tangent vector of this curve can be evaluated using the chain rule 
\begin{equation}
\overset{.}{\gamma }=\cosh \varphi _{2}N_{t}+\cos \theta _{2}\sinh \varphi
_{2}N_{y}+\sinh \varphi _{2}\sin \theta _{2}N_{z}  \tag{3.6}
\end{equation}%
and its magnitude $V_{2}$ is the speed, which is just the time rate of
change of the arc length along the curve $\gamma $. Hence, by using theorem
6 and theorem 7, $V_{2}^{y^{\ast }}=f_{1}V_{2}^{y}$ $=V_{2}\cos \theta
_{2}\sinh \varphi _{2}$ and $V_{2}^{z^{\ast }}=f_{2}V_{2}^{z}$ $=V_{2}\sinh
\varphi _{2}\sin \theta _{2}$ are just the radial vertical velocity while $%
V_{2}^{t}$ is the horizontal angular velocity and $V_{2}^{t^{\ast
}}=V_{2}^{t}$ $=V_{2}\cosh \varphi _{2}$ is the horizontal component of the
velocity vector. Similarly, $V_{2}^{y^{\ast }}$ is along the first axis, $%
V_{2}^{z^{\ast }}$ is along the second axis and $V_{2}^{t^{\ast }}$ is along
the third axis. Therefore, the specific kinetic energy can be given as 
\begin{equation*}
E_{_{\substack{ specific \\ energy}}}^{2}=\frac{1}{2}V_{2}^{2}=\frac{1}{2}%
\left( -V_{2}^{2}\cos ^{2}\theta _{2}\sinh ^{2}\varphi _{2}-V_{2}^{2}\sinh
^{2}\varphi _{2}\sin ^{2}\theta _{2}+V_{2}^{2}\cosh ^{2}\varphi _{2}\right) 
\end{equation*}%
\begin{equation}
=\frac{1}{2}\left( f_{1}\frac{dy}{ds}\right) ^{2}+\frac{1}{2}\left( f_{2}%
\frac{dz}{ds}\right) ^{2}-\frac{1}{2}\left( \frac{dt}{ds}\right) ^{2}, 
\tag{3.7}
\end{equation}%
by using the right side of the previous equations, the specific energy and
speed are constant along a geodesic. That is, its energy and hence specific
kinetic energy $E^{2}$ are constant and the speed $V_{2}=\sqrt{2E^{2}}$ is
constant along a geodesic, the velocities $V_{2}^{y^{\ast }}$ $=V_{2}\cos
\theta _{2}\sinh \varphi _{2}=f_{1}\frac{dy}{ds}$, $V_{2}^{z^{\ast }}$ $%
=V_{2}\sinh \varphi _{2}\sin \theta _{2}=f_{2}\frac{dz}{ds}$ and $%
V_{2}^{t^{\ast }}=V_{2}\cosh \varphi _{2}$ $=\frac{dt}{ds}$ are in the
angular direction multiplied by the radius $f_{2}$ and $f_{1}$ and from the
second geodesic equation, one writes 
\begin{equation*}
l_{_{\substack{ specific\text{ }angular \\ momentum}}}=\frac{\partial
L_{2}^{2}}{\partial \overset{.}{t}}=-2\overset{.}{t}=-2\cosh \varphi
_{2}V_{1}=-2\cosh \varphi _{2}\sqrt{2E^{2}}\Rightarrow \frac{-l_{2}}{2}=%
\overset{.}{t}, 
\end{equation*}
one can write the angular velocities $dy/ds$ and $dz/ds$ in the specific energy
formula according to the constant specific angular momentum and the radial motion and another constant of the motion is obtained as follows 
\begin{equation}
E_{_{\substack{ specific  \\ energy}}}=\frac{V_{2}^{2}}{2}\left( \sinh
^{2}\varphi_{2} -\frac{l_{2}^{2}}{8}\right).  \tag{3.8}
\end{equation}

3) For the elliptic surface of rotation $\Upsilon ^{3}$; if one wants to
obtain specific energy equations on this surface, one has to think the
integral of the length of the tangent vector of the
curve $\Upsilon ^{3}(\beta \left( s\right) ,\theta \left( s\right) ,t(s))$,
then the speed function is given as follows 
\begin{equation*}
I_{1}^{3}=\int ds=\int \frac{ds}{d\pi }d\pi =\int \sqrt{\left( \frac{d\beta 
}{d\pi }\right) ^{2}+\left( \frac{d\theta }{d\pi }\right) ^{2}+\left( \frac{%
dt}{d\pi }\right) ^{2}}d\pi 
\end{equation*}%
and this can be write as integral of half the length squared of the tangent vector, one gets 
\begin{equation*}
I_{2}^{3}=\frac{1}{2}\int \left( \frac{ds}{d\pi }\right) ^{2}d\pi =\frac{1}{2%
}\int \left( -\left( f_{2}\frac{d\beta }{d\pi }\right) ^{2}+\left( f_{4}%
\frac{d\theta }{d\pi }\right) ^{2}-\left( \frac{dt}{d\pi }\right)
^{2}\right) d\pi ,
\end{equation*}%
and the second Lagrangian function is called as  the energy function and is written as 
\begin{equation*}
L_{2}^{3}=\left( \beta ,\theta ,t,\frac{d\beta }{d\pi },\frac{d\theta }{d\pi 
},\frac{dt}{d\pi }\right) =-\frac{1}{2}\left( f_{2}\overset{.}{\beta }%
\right) ^{2}+\frac{1}{2}\left( f_{4}\overset{.}{\theta }\right) ^{2}-\frac{1%
}{2}\overset{.}{t}^{2}=E^{3},
\end{equation*}%
since the first Lagrangian $L_{1}^{3}$ is speed function one can write as 
\begin{equation*}
L_{1}^{3}=\left( \beta ,\upsilon ,t,\overset{.}{\beta },\overset{.}{\upsilon 
},\overset{.}{t}\right) =\sqrt{-\left( f_{2}\overset{.}{\beta }\right)
^{2}+\left( f_{4}\overset{.}{\theta }\right) ^{2}-\overset{.}{t}^{2}},
\end{equation*}%
and with the second Lagrangian, the angular equation is directly given the
constancy of the angular momentum $l_{3}$. Also, to derivative of tangent vector along $\Upsilon ^{3}(\beta
\left( s\right) ,\theta \left( s\right) ,t(s))$, by using the product and
chain rules, the tangent vector is obtain as  
\begin{equation}
\overset{.}{\gamma }=\cos \varphi _{3}N_{t}+\sin \varphi _{3}\cosh \theta
_{3}N_{\beta }+\sinh \theta _{3}\sin \varphi _{3}N_{\theta }.  \tag{3.9}
\end{equation}

Also, the tangent vector(velocity) of the geodesic on $\Upsilon ^{3}$ is
written as 
\begin{equation*}
\overset{\rightarrow }{V_{3}}=\frac{d\Upsilon ^{3}}{ds}=V_{3}^{\beta
}\Upsilon _{\beta }^{3}+V_{3}^{\theta }\Upsilon _{\theta
}^{3}+V_{3}^{t}\Upsilon _{t}^{3}
\end{equation*}%
and its magnitude $V_{3}$ is the speed. Also, by using theorem 8 and theorem
9, $V_{3}^{\beta ^{\ast }}=f_{2}V_{3}^{\beta }$ $=V_{3}\sin \varphi
_{3}\cosh \theta _{3}$ and $V_{3}^{\theta ^{\ast }}=f_{4}V_{3}^{\theta }$ $%
=V_{3}\sinh \theta _{3}\sin \varphi _{3}$ are the radial velocity while $%
V_{3}^{t}$ is the horizontal angular velocity. Then $V_{3}^{t^{\ast
}}=V_{3}^{t}$ $=V_{3}\cos \varphi _{3}$ is the horizontal component of the
velocity vector. Here, $V_{3}^{\beta ^{\ast }}$ is written along the first axis, $%
V_{3}^{\theta ^{\ast }}$ is written along the second axis and $V_{3}^{t^{\ast }}$ is
along the third axis.

Similarly, the angles $\theta _{3}$ and $\varphi _{3}$ give the direction of
the velocity according to the direction $\Upsilon _{\beta ^{\ast }}^{3}$.
Also, the speed is constant along the geodesic. Therefore, the specific
kinetic energy can be written as follows 
\begin{equation*}
E_{_{\substack{ specific \\ energy}}}^{3}=\frac{1}{2}V_{3}^{2}=\frac{1}{2}%
\left( V_{3}^{2}\sin ^{2}\varphi _{3}\cosh ^{2}\theta _{3}-V_{3}^{2}\sinh
^{2}\theta _{3}\sin ^{2}\varphi _{3}+V_{3}^{2}\cos ^{2}\varphi _{3}\right) 
\end{equation*}%
\begin{equation}
=-\frac{1}{2}\left( f_{2}\frac{d\beta }{ds}\right) ^{2}+\frac{1}{2}\left(
f_{4}\frac{d\theta }{ds}\right) ^{2}-\frac{1}{2}\left( \frac{dt}{ds}\right)
^{2}.  \tag{3.10}
\end{equation}%

Physically, the specific energy of the particle is constant because of its motion in space. Since its specific kinetic energy $E^{3}$ is
constant and the speed $V_{3}=\sqrt{2E^{3}}$ is constant along a geodesic.
Hence, $V_{3}^{\beta ^{\ast }}$ $=V_{3}\sin \varphi _{3}\cosh \theta
_{3}=f_{2}\frac{d\beta }{ds}$, $V^{\theta ^{\ast }}$ $=-V_{3}\sinh \theta
_{3}\sin \varphi _{3}=f_{4}\frac{d\theta }{ds}$ and $V_{3}^{t^{\ast
}}=V_{3}\cos \varphi _{3}$ $=\frac{dt}{ds}$ are velocities in the angular
direction multiplied by the radius $f_{2}$ and $f_{4}.$ Physically, by
thinking the second geodesic equation given as 
\begin{equation*}
l_{3}=\frac{\partial L_{2}^{3}}{\partial \overset{.}{t}}=-2\overset{.}{t}%
=-2\cos \varphi _{3}V_{3}=-2\cos \varphi _{3}\sqrt{2E^{3}}\Rightarrow \frac{%
-l_{3}}{2}=\overset{.}{t},
\end{equation*}%
and by using the variable angular velocities $d\beta /ds$, $d\theta /ds$ and for the radial motion and another constant of the motion the specific energy formula are written as 
\begin{equation}
E_{_{\substack{ specific \\ energy}}}=-\frac{V_{3}^{2}\sin ^{2}\varphi _{3}}{%
2}-\frac{l_{3}^{2}}{8}.  \tag{3.11}
\end{equation}

\begin{figure}[htp]
\centering
\subfloat[ \label{fig14e}] {\
\includegraphics[width=7cm,height=4.5cm]{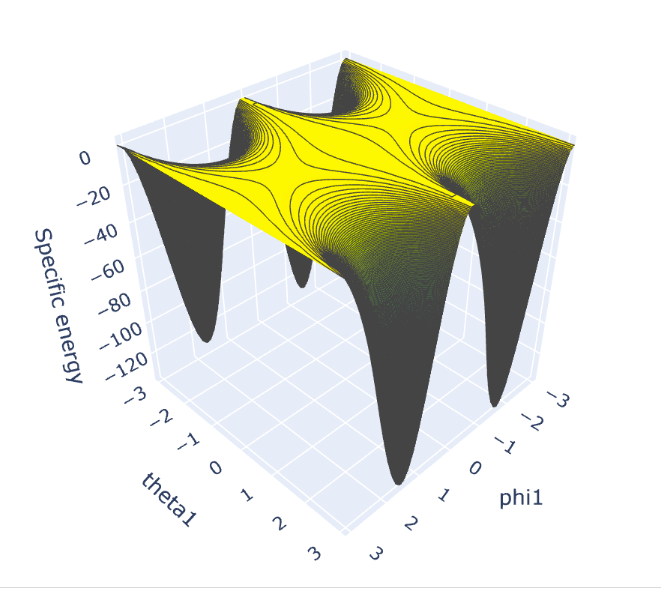} } \hspace*{.1cm}
\newline\caption{{The specific energy on
hyperbolic rotational surface $ \Upsilon^{1} $generated by the curve $\gamma(s)=(sins,0,0,coss)$} }%
\subfloat[  \label{figS23e2}] {\
\includegraphics[width=10cm,height=4.5cm]{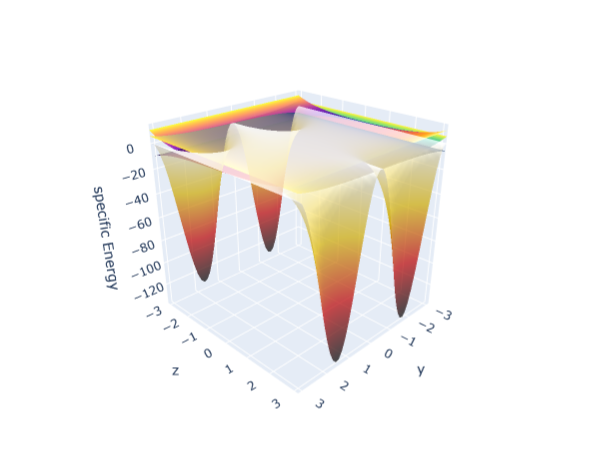}} \newline\caption{{The specific energy on hyperbolic rotational surface $ \Upsilon^{2} $ generated by the curve $\gamma(s)=(sins,coss,0,0)$} }%
\label{Fig1}%
\subfloat[ \label{figS56e2}] {\
\includegraphics[width=10cm,height=5.5cm]{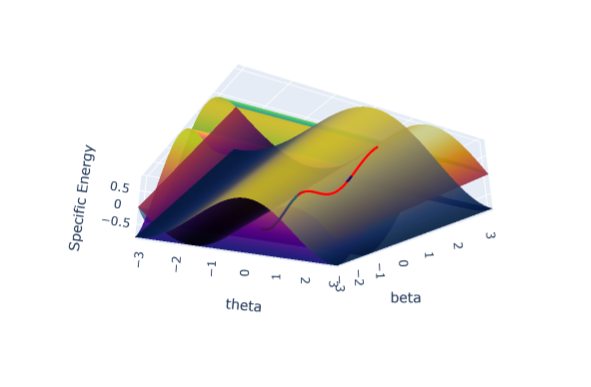} } \hspace*{.1cm}
\newline\caption{{The specific energy on elliptic rotational surface $ \Upsilon^{3} $ generated by the curve $\gamma(s)=(0,coss,0,coss)$} }%
\label{Fig2}%
\end{figure}

\section{Conclusion}


In this paper, the specific energy and the specific angular momentum on the surfaces of rotation are expressed in $E_{2}^{4}$ using the conditions of being geodesic, in which the curves can be chosen to be time-like curves, which allows us to constitute the specific energy and specific angular momentum. 





\section*{Acknowledgments}

The author wishes to express their thanks to the authors of
literatures for the supplied scientific aspects and idea for this study. 

\section*{Funding}

Not applicable.


\section*{Conflict of interest}

The author declares no conflicts of interest associated with this
manuscript. 





\end{document}